\documentstyle[11pt]{article}
\begin{document}
{\noindent \small Dokl. Akad. Nauk  388 (2003), no. 4, 439--442.}

\large

\

\begin{center}
 \quad UDC 517.925.42 \hfill MATHEMATICS \quad  \ \\
%\vskip0.3cm
 \bf A NEW APPROACH TO THE THEORY \\ OF
   ORDINARY DIFFERENTIAL EQUATIONS \\ WITH SMALL PARAMETER
 \footnote{Supported by RFBR grants 02-01-00189 and 02-01-00307}\\
\end{center}
\centerline{ \bf  M. I. Kamenskii, O. Yu. Makarenkov, P. Nistri}
\vskip0.5cm

In this paper we study both the periodic problem and the Cauchy problem
associated to the system of ordinary differential equations described by
\begin{equation} \label{1}
  \dot x=\varepsilon\phi(t,x)+\psi(t,x),
\end{equation}
where $\phi,\ \psi:{\rm R}\times{\rm R}^k\to {\rm R}^k$
are continuously differentiable $T$-periodic, with
respect to time $t$, functions and $\varepsilon$ is a small positive parameter.
Systems of the form (\ref{1}) represent a classic topic of
the theory of differential equations depending on a small parameter and such systems
have been investigated by a number of different methods. Here we devote the attention
only to the methods and results which are directly related to the aim of the present work.
In the papers \cite{Bobilev}, \cite{Krasn3}, \cite{KrasnZabr}
this problem is tackle by means of the theory of the rotation number of vector fields or
by theorems based on Poincar\`e index.
In this case for the $T$-periodic problem one usually assumes that
the system (\ref{1}) has isolated $T$-periodic solution $x_T$
of non-zero topological index at $\varepsilon=0$ , while for the Cauchy problem one assumes
the uniqueness of the solution $x_0$ defined on $[0,d].$ This second condition ensures that
the topological index of $x_0$ is also non-zero (see \cite{Perov}).
After this, by using results concerning the continuous dependence
of solutions on the operator equations on $\varepsilon$
(see, for instance, \cite{KrasnZabr}),
it is shown the existence and convergence of the
solutions to (\ref{1}) when $\varepsilon\ge 0$ sufficiently small to $x_T$ and $x_0$ respectively.
Observe, that for the Cauchy problem the closeness of the solutions is proved on the
interval where it is assumed the existence and uniqueness of the
solution $x_0.$
An other approach to deal with the periodic and Cauchy problem for system (\ref{1}) is the
averaging method proposed by N. N. Bogolubov - N. M. Krilov (see \cite{Bogolubov}).
To use this method system (\ref{1}) is first reduced to the standard form
\begin{equation} \label{Normal}
  \dot x =\varepsilon f(t,x),
\end{equation}
where function $f$ is $T$-periodic with respect to the first variable.
Then, topological methods and vector field theory can be applied to investigate system (\ref{Normal}),
(see, for instance,
\cite{AhmKam}, \cite{Kamenskii}, \cite{Samoylenko},
\cite{Strigin}).
For this, the following auxiliary system is considered
\begin{equation} \label{Usr}
  \dot x =\varepsilon f_0(x),
\end{equation}
where
$$
f_0(\xi)=\frac{1}{T}\int_0^T f(s,\xi)d\tau, \quad \xi\in {\rm R}^k.
$$
The isolated equilibrium states of system (\ref{Usr}) which have non-zero
topological index with respect to the vector field $-f_0$  give rise to $T$-periodic
solutions of system (\ref{Normal}). While the solutions of the Cauchy
problem for system (\ref{Normal}) are close to the corresponding
solutions of the Cauchy problem for system (\ref{Usr})
on the interval of length $(1/ \varepsilon) d$.
The main method for reducing system (\ref{1}) to standard form
(\ref{Normal}) consists in the following change of variable
\begin{equation} \label{Zamena}
z(t)=\Omega(0,t,x(t)),
\end{equation}
where $\Omega(\cdot,t_0,\xi)$ denotes the solution $x$
of (\ref{1}) with the initial condition $x(t_0)=\xi$ and $\varepsilon=0.$
Therefore, it is necessary to assume that the change of variable (\ref{Zamena}) is $T$-periodic
with respect to $t$ for every $T$-periodic function $x$
in order to use the classical averaging principle
(see \cite{Schneider}, \cite{Boston}).

In the present paper we do not require such $T$-periodicity condition for (\ref{Zamena}).
In fact, we assume only that the boundary of some open set $U\subset\rm R^k$ when $\varepsilon = 0$
represents the initial values of $T$-periodic solutions of system (\ref{1}).
This assumption implies that the $T$-periodic solutions of system (\ref{1}) at $\varepsilon=0$
are not isolated and so it is impossible to use directly results concerning
the continuous dependence of the solutions of the operator equations on $\varepsilon$.
This situation takes place, for instance, if system (\ref{1}) for  $\varepsilon=0$
is autonomous and it has an isolated cycle $x_0.$ In this case, as set $U$ can be
taken the interior of $x_0$. Such systems are very important in the
applications, but the investigation of the existence of $T$-periodic solutions
is difficult since the topological index of the set of
$T$-periodic solutions which arise from
cycle $x_0$ is equal to zero (see \cite{Bobilev7}).

In order to investigate the behaviour of solutions of system (\ref{1}) when
$\varepsilon\ge 0$ is small we introduce the following linear system
\begin{equation} \label{2}
\dot y=\phi(t,\Omega(t,0,\xi))+\psi'_{(2)}(t,\Omega(t,0,\xi))y,
\end{equation}
where $\xi\in {\rm R}^k$.
We will establish our results in terms of system (\ref{2}). The classical
results of the averaging theory both of the $T$-periodic problem and the
Cauchy problem will follow from our results. We will also show that the result
from  \cite{Koddington} in the case when $\psi(t,x)=Ax$, where the
matrix $A$ has two simple eigenvalue $\pm i\frac{T}{2\pi}$, is a
consequence of our results.

\vskip0.5truecm

1. We consider first the problem of the existence of $T$-periodic solutions.
We denote by $\eta(\cdot,s,\xi)$ the solution $y$ of system (\ref{2})
satisfing $y(s)=0.$ In this section we will use the classical notion of the
rotation number of a continuous map $F:\bar U\to \bar U$ defined in the closure
of a bounded set $U\subset {\rm R}^k$ (see \cite{KrasnZabr}) and we denote it by
$\gamma (F,U)$.

\vskip0.5truecm
\noindent
{\bf Theorem 1.}

Let the set $U\subset {\rm R}^k$ be open and bounded. Assume, that the
following conditions hold
\begin{eqnarray}
& {\rm (A0)}&\quad \Omega(T,0,\xi)=\xi, \quad \xi\in \partial U, \nonumber
\\
&{\rm (A1)}&\quad \eta(T,s,\xi)-\eta(0,s,\xi)\not=0,\quad s\in[0,T],\
\xi\in\partial U, \nonumber \\
&{\rm (A2)}& \quad \gamma(\eta(T,0,\cdot),U)\not=0. \nonumber
\end{eqnarray}
Then there exists $\varepsilon _0\ge 0$ such, that
system (\ref{1}) has at least one $T$-periodic solution belonging to the set
$X=\left\{x:\Omega(0,t,x(t))\in U,\ t\in[0,T]\right\}$ for all
$\varepsilon\in(0,\varepsilon_0).$

\vskip0.5truecm

We would like to point out that there are many cases for which it is enough to verify condition (A2)
for some particular function $\phi,$ which does not necessarly coincide with the one given in (\ref{1}).

\vskip0.5truecm
\noindent
{\bf Theorem 2.} Assume condition (A0).
Then for every functions $\phi_1$ and $\phi_2,$ such that
the function
$\phi_{\lambda}=\lambda\phi_1+(1-\lambda)\phi_2$
satisfies condition (A1) for all $\lambda\in[0,1]$, the rotations
$\gamma(\eta_1(T,0,\cdot),U)$ and $\gamma(\eta_2(T,0,\cdot),U)$
coincide.

\vskip0.5truecm

Assume that the dimension of the phase space $k$ is equal to $2$
and the system (\ref{1}) has an  $T$-periodic cycle, that is
it has $T$-periodic solution $x_0,$  such that function $x_0(t+\theta)$
on $t$
is a solution of (\ref{1}) for $\varepsilon=0$ and $\theta\in[0,T].$
Assume that the cycle $x_0$ is simple, that is curve
$x_0$ has no self-joints. Than by Djordan theorem
the curve $x_0$ restricts an one-connected domain $U$ of the
space ${\rm R}^2.$
By using the fact that the rotation number of the field of tangents and
the rotation number of the field of normals are equal to 1 on the cycle
we derive in the sequel an existence result for $T$-periodic solutions.
Since domain $U$ is one-connected than
it homeomorphic to the unique circle $B_1.$
Denote by $g$ some homeomorphism of $U$ on $B_1.$ For an arbitrary
$\delta\in[0,1)$ define $\delta$-contraction $W_{\delta}(U)$ of the
domain $U$ by the formula $W_{\delta}(U)=g^{-1}((1-\delta)g(U)).$ Now if

(A3) $T$-periodic system \ $\dot y=\psi'_{(2)}(t,x_0(t+\theta))y$ \ has 1 as a simple

$\qquad$ Floquet exponent for all $\theta\in[0,T]$
\\
than the points of boundary of the set $W_{-\delta}(U)$ are
the points of $T$-irreversibility (see \cite{Krasn3}) of the
system (\ref{1}) for $\varepsilon=0.$ Thus by condition (A3) the system
(\ref{1}) has $T$-periodic solution acting to $W_{-\delta}(U)$ for
small $\varepsilon>0$ and $\delta>0$ (see \cite{Krasn3}, theorem 6.1).
The theorem 6.1 from \cite{Krasn3} is not applicable when the condition (A3)
does not hold generally speaking.
In this case a condition of existence of $T$-periodic solutions
for system (\ref{1}) can be derived by the ground of theorems 1 and 2.
In fact, observing than in this case
$<\eta(T,s,x_0(\theta))-\eta(0,s,x_0(\theta)),{\dot x_0(\theta)}^{\bot}>\not=0$
for such and only such $s,\theta\in[0,T]$ that
\begin{eqnarray}
  {\rm (A3)}_1 & & \int_0^T {\rm e}^{-\int_0^t {\rm Sp}\psi'_{(2)}(\tau,x_0(\tau))d\tau}<\phi(t-\theta,x_0(t)),{\dot x_0(t)}^{\bot}>dt\not=0, \nonumber
\end{eqnarray}
where $<\cdot,\cdot>$ is a scalar product in ${\rm R}^2$ and ${\dot x_0(\theta)}^{\bot}$
is a rotation of ${\rm R}^2$ and ${\dot x_0(\theta)}^{\bot}$ against clockwise
we can derive the following establishment.

\vskip0.5truecm

{\bf Theorem 3.} Let $x_0$ is a simple cycle of the system (\ref{1}) and $U$
is it's interior. Assume the condition ${(A3)}_1$ holds for all
$s,\theta\in[0,T].$ Then for sufficiently small $\varepsilon>0$ the system
(\ref{1}) has $T$-periodic solution acting in $U.$

\vskip0.5truecm

Theorem 3 can be generalized to the case of systems (\ref{1}) of
even dimension $k=2p,$ namely to the case of systems (\ref{1})
consisting of $p$ two-dimensional systems which have simple
cycles $x_1$,...,$x_p$ of the same period $T.$ In this case one might
take as $U$ the cartesian product of the interiors of this cycles in the corresponding spaces.

\vskip0.5truecm

2. Let us now to study the Cauchy problem. In this section the initial
condition for system (\ref{1}) is fixed. Assume, that the limit
\begin{equation} \label{sred}
\Phi(\xi)=-\lim_{n\to\infty}\frac{\eta(-nT,0,\xi)}{nT}
\end{equation}
exists uniformly with respect to $\xi\in B(0,r)$ for every $r>0.$
The limit (\ref{sred}) exists, for instance, in the case when the system
(\ref{2}) is stable for $t\to -\infty$ uniformly with respect to
$\xi\in B(0,r)$ for every $r> 0.$

\vskip0.5truecm
\noindent
{\bf Theorem 4.}
Assume that the system $\dot z=\Phi(z)$ has an unique solution $z_0$ on the interval
$[0,d]$, satisfing $z(0)={\xi}_0.$ Then for every $\gamma > 0$
there exists $\varepsilon_0 > 0$ such that (\ref{1}) has a
solution satisfing $x_{\varepsilon}(0)={\xi}_0$ and such that
$$
  \left\|x_{\varepsilon}(t)-\Omega(t,0,z(\varepsilon t))\right\|\le
\gamma,\quad t\in\left[0,\frac{d}{\varepsilon}\right]
$$
for every $\varepsilon\in (0,\varepsilon_0).$

\vskip0.5truecm

Classical results of averaging theory given by N. N. Bogolubov - N. M. Krilov
for $T$-periodic problem and for Cauchy problem (see \cite{Bogolubov})
follow from our results.
To see this, it is sufficient to put $\psi=0$ in (\ref{1}) and use the formula
$$
  \eta(\pm nT,s,\xi)-\eta(0,s,\xi)=\int_{s\mp nT}^s \phi(\tau,\xi)d\tau=
  \int_{\mp nT}^0 \phi(\tau,\xi)d\tau.
$$
in order to compare the conditions of theorems 1 and 4 with the classical ones.

We now describe how from theorem 1 we can obtain a classical result
from (\cite{Koddington}, Theorem 3.1, p. 362) in the case when
$\psi(t,x)=Ax=(-x_2,x_1)$ and $\phi(t,x)=(0,g(t,-x_1,x_2))$, where the
function $g$ is $2\pi$-periodic with respect to the first variable.
Thus the dimension of the phase space is $2$ and the matrix $A$ has simple
eigenvalues $\pm i$. In this example we have choosen the notations in such
a way that formulas from \cite{Koddington} and formulas from theorem 1
coincide. We have
$$
\eta(T,s,\xi(a,\theta))-\eta(0,s,\xi(a,\theta))=
\left( \begin{array}{cc}
         \cos\theta \sin\theta \\
         -\sin\theta \cos\theta \end{array} \right)H(a,\theta),
$$
where $\xi(a,\theta)=(-a\cos\theta,a\sin\theta)$ and
$$
 H(a,\theta)=\int_0^{2\pi}\left(\begin{array}{c}
    (\sin\tau)f(\tau+\theta,a\cos\tau,-a\sin\tau) \\
    (\cos\tau)f(\tau+\theta,a\cos\tau,-a\sin\tau)
   \end{array} \right) d\tau.
$$

It is assumed in \cite{Koddington} the existence
of $a_0,{\theta}_0\in {\rm R}$ such that \\
\centerline{$H(a_0,{\theta}_0)=0$ \ and \ ${\rm det} | H'(a_0,{\theta}_0)|\not=0.$}\\
But by this assumption there exists such the set
$
V=\left(a_0-\delta,a_0+\delta\right) \times \times
   \left({\theta}_0-\delta,{\theta}_0+\delta\right),$
that $|\gamma(",V)|=1$.
Since vector fields $\eta(T,s,\xi(\cdot,\cdot))-\eta(0,s,\xi(\cdot,\cdot))$ and
$H(\cdot,\cdot)$ are homotopic on $\partial V$ for every $s\in[0,2\pi]$ then
these fields have the same rotation number on $V.$
Without loss of generality we can consider $0\le\delta \le a_0$ and $0\le\delta\le \pi.$
In this case the continuous function $\xi$ maps the open set $V$ on the open
set $\xi(V).$ Hence $\gamma(\eta(T,0,\cdot)-\eta(0,0,\cdot),\xi (V) )=
\gamma(\eta(T,s,\xi(\cdot,\cdot))-\eta(0,s,\xi(\cdot,\cdot)),V)$
and by theorem 1 system (\ref{1}) has $2\pi$-periodic solutions in the
set $X=\left\{x:(a,t+\theta)\in V,\ (a,\theta)=\xi^{-1}(x(t)),\
t\in[0,2\pi]\right\}$ for sufficiently small $\varepsilon\ge 0.$

\end{document}